\documentclass[a4paper,12pt]{article}

\usepackage[english]{babel}
\usepackage{mathrsfs}
\usepackage{amsfonts}
\usepackage[centertags]{amsmath}
\usepackage{amsthm}
\usepackage{enumitem}
\usepackage{latexsym,amsmath,amssymb}
\usepackage{graphicx,color}
\usepackage[all]{xy}
\newcommand{\ds}{\displaystyle}
\newcommand{\cua}{^{2}}
\newcommand{\bra}{\langle}
\newcommand{\ket}{\rangle}

\newcommand{\R}{\mathbb{R}}      

\setenumerate[1]{label=(\roman*), ref=(\roman*)}
\usepackage[all]{xy}

\newtheorem{thm}{Theorem}[section]
\newtheorem{prop}[thm]{Proposition}

\newtheorem{lemma}[thm]{Lemma}

\newtheoremstyle{obs}
  {3pt}
  {3pt}
  {}
  {}
  {\bfseries}
  {.}
  {.5em}
  {}
\theoremstyle{obs}
\newtheorem{remark}[thm]{Remark}
\newtheorem{example}[thm]{Example}
\newtheorem{defn}[thm]{Definition}

\newcommand{\todo}[1]{\vspace{5 mm}\par \noindent
\marginpar{\textsc{ToDo}}
\framebox{\begin{minipage}[c]{0.95 \textwidth}
\tt #1 \end{minipage}}\vspace{5 mm}\par}
\newcommand{\lp}{\left(}
\newcommand{\rp}{\right)}

\def\qed{\ifvmode\removelastskip\fi
{\unskip\nobreak\hfil\penalty50\hbox{}\nobreak\hfil \hbox{\vrule
height1.2ex width1.2ex}\parfillskip=0pt \finalhyphendemerits=0
\par \smallskip}}

\title{Morse families in optimal control problems}
\author{\textsc{Mar\'ia Barbero-Li\~n\'an}\thanks{mbarbero@math.uc3m.es}\\
\small
Departamento de Matem\'aticas, Universidad Carlos III de Madrid, \\ \small Avenida de la Universidad 30, 28911 Legan\'es, Madrid, Spain
\\ \small and Instituto de Ciencias Matem\'aticas (CSIC-UAM-UC3M-UCM)\\ \and 
\textsc{David Iglesias Ponte}\thanks{diglesia@ull.es} \\
\small Unidad asociada ULL-CSIC ``Geometr\'{\i}a Diferencial y Mec\'anica Geom\'etrica" \\
\small Departamento de Matem\'atica Fundamental, Facultad de Matem\'aticas, \\ \small Universidad de la Laguna, Spain
\\
 \and
\textsc{David Mart\'{\i}n de Diego}\thanks{david.martin@icmat.es} \\
\small
Instituto de Ciencias Matem\'aticas (CSIC-UAM-UC3M-UCM) \\ \small  C/Nicol\'as
Cabrera 13-15, 28049 Madrid, Spain
}

\begin{document}

\maketitle

\begin{abstract}
We geometrically describe optimal control problems in terms of Morse families in the Hamiltonian framework.
These geometric structures allow us to recover the classical first order necessary conditions for optimality and the starting point to run an integrability algorithm. Moreover the integrability algorithm is adapted to optimal control problems in such a way that the trajectories originated by discontinuous controls are also obtained. From the Hamiltonian viewpoint we obtain the equations of motion for optimal control problems in the Lagrangian formalism by means of a proper Lagrangian submanifold. Singular optimal control problems and overdetermined ones are also studied along the paper.
\end{abstract}

\section{Introduction}

In the late fifties the interest in optimal control problems grew amazingly due to the applications in military missions and later on
due to the applications in space missions. Pontryagin and his collaborators provided us with a useful result to characterize candidates to be optimal solutions~\cite{Pontryagin}. Since then great efforts have been made to understand optimal control problems from a differential geometric viewpoint (see  ~\cite{Agrachev,2005BulloLewis,1997Jurdjevic, 1998Sussmann} and references therein). It turned out that both symplectic~\cite{S98Free} and presymplectic~\cite{BarMu} formalisms were useful in that sense.

In this paper we present a novel formulation of optimal control problems by using different geometric tools (see~\cite{2012FerLeMa} and~\cite{GraGra2008} for a different but related approach in the scope of calculus of variations). Those tools are  Morse families~\cite{LiMarle}  and Lagrangian submanifolds~\cite{Weinstein}, which are connected by the fact that any Lagrangian submanifold is locally described by a Morse family.  We will show how these geometric elements, once properly adapted to optimal control theory, summarize a weaker version of the first order necessary conditions given by Pontryagin's Maximum Principle. At this stage a suitable adaptation of the integrability algorithm in~\cite{1995MMT} can be considered to find the set of candidates to be optimal solutions. Both Morse families in Section~\ref{Sec:Morse}  and Lagrangian submanifolds in Section~\ref{Sec:LagranOCP} describe a new setting to distinguish and discuss different kinds of optimal control problems. All these results are applied to different singular and real examples, among them overdetermined control systems. These last systems are very useful for some applications in order to gain mobility in the control system under study.

Moreover, it is important to stress that in our formulation the controls will play the role of parameters, and the spaces where  the associated optimal control equations are given are exactly the same ones as for classical Hamiltonian and Lagrangian  mechanics. In this sense our point of view of a control problem is similar to \cite{1997Jurdjevic}.

The outline of the paper is as follows. Section~\ref{Sec:Preliminaries} summarizes all the geometric tools and results necessary for the paper. Section~\ref{Sec:IntroOCP} briefly recalls the notion of optimal control problems and Pontryagin's Maximum Principle. Some of the main contributions of this paper are in Sections~\ref{Sec:Morse} and~\ref{Sec:LagranOCP}  where Morse families and Lagrangian submanifolds, respectively, are used to describe optimal control problems and their regularity. Several examples are provided to make the theory clearer. Another important constribution appears in Section~\ref{Sec:Integrab} where the integrability algorithm in~\cite{1995MMT} is adapted to optimal control problems with the purpose of recovering even the solutions originated by bang- bang controls.

\section{Geometric preliminaries}\label{Sec:Preliminaries}

In this section we briefly introduce all the definitions and results from differential geometry that are necessary in the sequel. See~\cite{AbMa,GuiStern,LiMarle,TuHamilton,Tu,Weinstein} for more details.

\subsection{Lagrangian submanifolds and Morse families}\label{SubSec:LagranMorseIntro}

We define the two main geometric tools that are used in this paper: Lagrangian submanifolds and Morse families. More details can be found in~\cite{AbMa,1989LeRo,GuiStern,LiMarle,Weinstein} and references therein.

Let us recall that a symplectic vector space is a pair $(E,\Omega)$ where $E$ is a vector space and $\Omega\colon E\times E \rightarrow \mathbb{R}$ is a skew-symmetric bilinear map of maximal rank. We distinguish the following type of vector  subspaces of a symplectic vector space:
\begin{defn} Let $(E,\Omega)$ be a symplectic vector space and
$F\subset E$ be a subspace. The $\Omega$-\textit{orthogonal complement
of $F$} is the subspace defined by
\begin{equation*}
F^\perp=\{e\in E \; | \; \Omega(e,e')=0 \; \mbox{for all } e'\in
F\}.
\end{equation*}
A subspace $F$ of a symplectic vector space is called
\begin{enumerate}
\item  \textit{isotropic} if $F\subset F^\perp$, that is,
$\Omega(e,e')=0$ for all $e,e'\in F$.
\item \textit{Lagrangian} if $F$ is isotropic and has an
isotropic complement, that is, $E=F\oplus F'$, where $F'$ is
isotropic.
\end{enumerate} \label{defn:lagrangianSubspace}
\end{defn}

A well-known characterization of Lagrangian subspaces of a finite dimensional symplectic vector space is the following one:

\begin{prop} Let $(E,\Omega)$ be a symplectic vector space and
$F\subset E$ a subspace. Then the following assertions are
equivalent:
\begin{enumerate}
\item $F$ is Lagrangian,
\item $F=F^\perp$,
\item $F$ is isotropic and ${\rm dim} \, F=\frac{1}{2}{\rm dim}\,
E$.
\end{enumerate}
\end{prop}

As a consequence, Lagrangian subspace $F$ of $E$ can be characterized by checking if it has half the dimension of $E$ and $\Omega_{|F}=0$.

 Remember that a symplectic manifold
$(M,\omega)$ is defined by a differentiable manifold $M$ and a non-degenerate closed 2-form $\omega$ on $M$. Therefore, for each $x\in M$, $(T_xM, \omega_x)$ is a symplectic vector space
and a symplectic manifold has even dimension. Denote by $\flat_\omega: TM\rightarrow T^*M$ the isomorphism $\flat_\omega (v)=i_v\omega$, with $v\in TM$, and by $\sharp_M= (\flat_\omega )^{-1}$.

The notion of Lagrangian subspace can be transferred to submanifolds by requesting that the tangent space of the
submanifold is a Lagrangian subspace for every point in the submanifold of a symplectic manifold.

\begin{defn} Let $(M,\omega)$ be a symplectic manifold and ${\rm
f}\colon N\rightarrow M$ an immersion. It is said that $N$ is a
\textit{Lagrangian immersed submanifold} of $(M,\omega)$ if $({\rm
T}_x{\rm f})({\rm T}_xN) \subset {\rm T}_{{\rm f}(x)}M$ is a Lagrangian subspace
for each $x\in N$. We say that ${\rm f}$ is a \textit{Lagrangian immersion}.
\label{defn:lagrangianSubm}
\end{defn}

Note that ${\rm f}\colon N\rightarrow M$ is isotropic if and only if
${\rm f}^*\omega=0$, that is, $\omega({\rm T}_x{\rm f} (u_x), {\rm T}_x{\rm f} (v_x))=0$ for every $u_x,v_x\in T_xN$ and for every $x\in N$.

The canonical model of symplectic manifold is the cotangent bundle $T^*Q$ of an arbitrary manifold $Q$.
Denote by $\pi_Q\colon T^*Q \rightarrow Q$ the canonical projection and define a canonical 1-form $\theta_Q$ on $T^* Q$ by
\begin{equation*}
 \left(\theta_Q\right)_{\alpha_q}(X_{\alpha_q})=\langle \alpha_q, {\rm T}_{\alpha_q} \pi_Q(X_{\alpha_q})\rangle,
\end{equation*}
where $X_{\alpha_q}\in T_{\alpha_q}T^*Q$, $\alpha_q\in T^*Q$ and $q\in Q$. If we consider
bundle coordinates $(q^i,p_i)$ on $T^* Q$ such that $\pi_Q(q^i,p_i)=q^i$, then
\begin{equation*}
\theta_Q=p_i  {\rm d}q^i\, .
 \end{equation*}
 The 2-form $\omega_Q=-{\rm d}\theta_Q$ is a symplectic form on $T^*Q$ with local expression
\begin{equation*}
\omega_Q={\rm d}q^i \wedge {\rm d}p_i.
\end{equation*}
The Darboux theorem states that this is the local model for an arbitrary symplectic manifold $(M,\omega)$: there exist local
coordinates $(q^i,p_i)$ in a neighbourhood of each point in $M$ such that $\omega={\rm d}q^i \wedge {\rm d}p_i$.

There are different ways to define Lagrangian submanifolds as we will show in the sequel.
A relevant example of Lagrangian submanifold of the cotangent bundle is the following one.
For complete proofs of the following statements we refer the reader to~\cite{LiMarle}.

\begin{prop} Let $\alpha$ be a one-form on $Q$. The submanifold ${\rm Im} \, \alpha$ is a Lagrangian submanifold  if and only if
$\alpha$ is closed. \label{prop:LagrangianClosedForm}
\end{prop}

Given  a symplectic manifold $(M, \omega)$, $\dim M=2n$, it is well-known that the tangent bundle $TM$ is equipped with a
symplectic structure denoted by $\mathrm{d}_T\omega$, where ${\rm d}_T \omega$ denotes the tangent lift of $\omega$ to $TM$. If we
take Darboux coordinates $(q^i,p_i)$ on $M$, that is, $\omega=\mathrm{d}q^i\wedge \mathrm{d}p_i$, then
$\mathrm{d}_T\omega=\mathrm{d}\dot{q}^i\wedge \mathrm{d}p_i+\mathrm{d}q^i\wedge \mathrm{d}\dot{p}_i$, where
$(q^i, p_i; \dot{q}^i, \dot{p}_i)$ are the induced coordinates on $TM$. Denoting the bundle coordinates on $T^*M$ by
$(q^i, p_i; a_i, b^i)$, then $\omega_{M}=\mathrm{d}q^i\wedge \mathrm{d}a_i+\mathrm{d}p_i\wedge \mathrm{d}b^i$.
If $\flat_\omega: TM\to T^*M$ is the isomorphism defined by $\omega$, that is, $\flat_{\omega}(v)={\rm i}_{v}\,\omega$, then $\flat_{\omega}(q^i, p_i; \dot{q}^i, \dot{p}_i)=(q^i, p_i; -\dot{p}_i, \dot{q}^i)$.

Given a function $H: M\to \R$, and its associated Hamiltonian vector field $X_H$, that is, ${\rm i}_{X_H}\omega=\mathrm{d}H$,
the image of  $X_H$, $\hbox{Im } X_H$, is a Lagrangian submanifold of  $(TM, \mathrm{d}_T\omega)$.
 Moreover, given a vector field $X$ on $M$, it is locally Hamiltonian if and only if its image $X(M)$ is a Lagrangian submanifold of $(TM, \mathrm{d}_T\omega)$.
It is interesting to note that
$\mathrm{d}_T\omega=-\flat_{\omega}^{*}\, \omega_M$ and $\flat_{\omega}(X_H(M))=\mathrm{d}H(M)$.

%
%

Another geometric element that can be used to define Lagrangian submanifolds is a Morse family. The notion of a Morse family or
phase function was introduced by L. H\"ormander~\cite{Hormander}.  Here we just give the key definitions and essential results for this paper using a similar notation to the one in \cite{LiMarle}.

\begin{defn} Let $\pi\colon B \rightarrow Q$ be a submersion of a differentiable manifold $B$ onto a differentiable manifold $Q$. Let $S\colon B\rightarrow \mathbb{R}$ be a differentiable function. The function $S$ is called a \textit{Morse family} if the image of the differential of $S$, $ {\rm d}S(B)\subset T^*B$, and the conormal bundle
\[
Q_\pi=(\ker \, {\rm T}\pi)^0=\left\{\alpha\in T_\mu ^*B\; |\; \langle \alpha, v\rangle=0, \hbox{ for all } v\in \ker T_{\pi_B(\mu)}\pi\right\}\subset T^*B
\]
  are transverse in $T^*B$, that is,
\begin{equation} \forall \alpha\in Q_\pi \cap  {\rm d}S(B)  \subseteq T^*B,   \quad {\rm T}_\alpha( {\rm d}S(B))+T_\alpha Q_\pi=T_\alpha (T^*B).
\label{Eq:Transverse}
\end{equation}
\label{Def:MorseFamily}
\end{defn}
Observe that under these transversality conditions $Q_\pi \cap  {\rm d}S(B)$ is an isotropic  submanifold of $(T^*B, \omega_B)$, since it is contained in the Lagrangian submanifold ${\rm d}S(B)$. Moreover
\[
\dim (Q_\pi \cap  {\rm d}S(B))=\dim Q  \, .
\]
If the submersion $\pi: B\rightarrow Q$ is expressed locally by $\pi(q^i, y^a)=(q^i)$, then $S: B\rightarrow \R$ is a Morse family in a local sense if the \[
\hbox{rank  of} \left(
\frac{\partial^2 S}{\partial q^i \partial y^b}, \frac{\partial^2 S}{\partial y^a \partial y^b}
\right) \hbox{ is maximal}
\]
for all $(q^i, y^a)$ satisfying $\displaystyle{\frac{\partial S}{\partial y^a}=0}$.

Define the  morphism $\mathfrak{j}\colon Q_\pi \rightarrow T^*Q$  by \begin{equation*} \langle \eta,v\rangle =\langle \mathfrak{j}(\eta), {\rm T}\pi (v) \rangle\end{equation*}
for all $\eta\in Q_\pi$ and for all $v\in T_{\pi (\eta)} B$. The following result will be useful in the sequel. For a proof see, for instance,~\cite[Appendix 7, Proposition 1.12]{LiMarle}.
\begin{prop}\label{prop:img:j:lagrangian}
Let $S\colon B\rightarrow \mathbb{R}$ be a Morse family. The restriction of the morphism $\mathfrak{j}\colon Q_\pi \rightarrow T^*Q$    to the isotropic submanifold $Q_\pi \cap {\rm d}S(B)$ is a Lagrangian immersion of $Q_\pi \cap {\rm d}S(B)$ in $(T^*Q, \omega_Q)$. This Lagrangian immersion is said to be generated by the Morse family $S$.
\label{Prop:MorseFamily}
\end{prop}

The assumptions in Proposition~\ref{Prop:MorseFamily} can be weakened by requiring that $Q_\pi$ and ${\rm d}S(B)$ are only weakly transverse, instead of the transversality required in the Definition~\ref{Def:MorseFamily} of a Morse family. Remember that two submanifolds are weakly transverse if $Q_\pi \cap {\rm d} S(B)$ is a submanifold of $T^*B$ and if, for all $x\in Q_\pi \cap {\rm d} S(B)$,
  \begin{equation}
  {\rm T}_x(Q_\pi \cap {\rm d} S(B))={\rm T}_x(Q_\pi)\cap {\rm T}_x ({\rm d}S(B)).\label{eq:CondWeakTrans}
  \end{equation}
 Thus, transversality implies weak tranversality.
Some authors use the name clean intersection instead of weak transversality.
These notions of transversality are defined for any pair of submanifolds, without being of a symplectic manifold.

 Assuming only clean intersection of $Q_\pi$ and ${\rm d}S(B)$, then  $\mathfrak{j}(Q_\pi \cap {\rm d} S(B))$ is not, in general, an immersed submanifold of $T^*Q$ 
 because it can admit  multiple points. These multiple points $x^*\in \mathfrak{j}(Q_\pi \cap {\rm d} S(B))$ are points at
which there exist several distinct vector subspaces of the space tangent to $T^*Q$ in such a way that each one is tangent to a $n$-dimensional
 submanifold contained in $\mathfrak{j}(Q_\pi \cap {\rm d} S(B))$. The following result is proved in ~\cite[Chapter III, Proposition 14.19]{LiMarle}.

  \begin{prop}
  Let $S\colon B\rightarrow \mathbb{R}$ be a Morse family.  Assume that $Q_\pi$ and ${\rm d}S(B)$ are weakly transverse. Then,  the restriction of $\mathfrak{j}$ to $Q_\pi \cap {\rm d} S(B)$ is of constant rank, and every point of $Q_\pi \cap {\rm d} S(B)$ has an open neighborhood $V$ in $Q_\pi \cap {\rm d} S(B)$ whose image $\mathfrak{j}(V)$ is a Lagrangian submanifold of $(T^*Q,\omega_Q)$. We say that $\mathfrak{j}(Q_\pi \cap {\rm d} S(B))$ is, in a generalized sense, an immersed Lagrangian submanifold of $(T^*Q,\omega_Q)$ which may include multiple points.  \end{prop}

The following proposition determines the relationship between Lagrangian submanifolds and Morse families as shown in
Remark~\ref{Rem:LagrangianMorse}.

\begin{prop}[\cite{GuiStern},\cite{LiMarle},\cite{Weinstein}] Let $\Sigma$ be a Lagrangian submanifold of $(T^*Q,-{\rm d} \theta_Q)$. We consider the following three properties:
\begin{enumerate}
\item The closed 1-form induced on $\Sigma$ by the Liouville 1-form $\theta_Q$ is exact.
\item There exists a Lagrangian subbundle $E$ of the symplectic bundle $(T_\Sigma(T^*Q),\omega_Q \left. \right|_\Sigma)$ (that is, the restriction to $\Sigma$ of the cotangent bundle $(T^*Q, \omega_Q)$)
which is a complement of both $T\Sigma$ and the restriction of the vertical tangent bundle of $T^*Q$ to $\Sigma$, denoted by $V_\Sigma(T^*Q)$.
\item There exists a surjective submersion $\pi\colon B \rightarrow Q$ and a Morse family $S\colon U \rightarrow \mathbb{R}$, the latter defined on an open subset $U$ of $B$, such that the Lagrangian immersion generated by $S$ is an embedding with image $\Sigma$.
\end{enumerate}
If properties (i) and (ii) are both satisfied, so is property (iii). Conversely, if property (iii) is satisfied, then property (i) also holds. \label{Prop:3Prop}
\end{prop}
\remark \label{Rem:LagrangianMorse} As a direct consequence we deduce that any Lagrangian submanifold is locally, i.e., in a neighborhood of each of
its points, the image of a Lagrangian immersion generated by a Morse family. See~\cite[Appendix 7, Proposition 9]{LiMarle} for more details.

\subsection{Tulczyjew triple}\label{Sec:Tulczy}

The theory of Lagrangian submanifolds gives an intrinsic geometric description of Lagrangian and Hamiltonian dynamics~\cite{TuHamilton},~\cite{Tu}. Moreover, it allows us to relate Lagrangian and Hamiltonian formalisms using as a main tool the so-called Tulczyjew's triple
\begin{equation*}
\xymatrix{ T^*TQ && TT^*Q \ar[rr]^{\hbox{\small{$\beta_Q$}}} \ar[ll]_{\hbox{\small{$\alpha_Q$}}} && T^*T^*Q\, .}
\end{equation*}

 The Tulczyjew map $\alpha_{Q}$ is an isomorphism between $TT^*Q$ and $T^*TQ$. Beside, it is also a symplectomorphism
between these  vector bundles considered as  symplectic manifolds, i.e. $(TT^*Q\,,\,\mathrm{d}_T\,\omega_Q)$, where
$\mathrm{d}_T\,\omega_Q$ is the tangent lift of $\omega_Q$, and $(T^{*}TQ, \omega_{TQ})$. For completeness, we recall
the construction of the symplectomorphism $\alpha_Q$. To do this, it is necessary to introduce the canonical involution $\kappa_Q$ on $TTQ$
\begin{equation*}
\xymatrix{
TTQ\ar[d]_{\tau_{TQ}}\ar[r]_{\kappa_{Q}}&TTQ\ar[d]^{T\tau_{Q}}\\
TQ\ar[r]_{\mbox{Id}}&TQ,
}
\end{equation*}
defined by
$$
\kappa_{Q}\lp\frac{{\rm d}}{{\rm d}s}\Big |_{s=0}\frac{{\rm d}}{{\rm d}t}\Big |_{t=0}\hspace{1mm}\chi\lp s,t\rp\rp=
\frac{{\rm d}}{{\rm d}s}\Big |_{s=0}\frac{{\rm d}}{{\rm d}t}\Big |_{t=0}\hspace{1mm}\tilde\chi\lp s,t\rp,
$$
where $\chi:\mathbb{R}\cua\rightarrow Q$ and $\tilde\chi:\mathbb{R}\cua\rightarrow Q$ are related by $\tilde\chi\lp s,t\rp=\chi\lp t,s\rp$. If $\lp q^{i}\rp$ are the local coordinates for $Q$, $\lp q^{i},v^{i}\rp$ for $TQ$ and $\lp q^{i},v^{i},\dot q^{i},\dot v^{i} \rp$ for $TTQ$, then the canonical involution is locally given by $\kappa_{Q}\lp q^{i},v^{i},\dot q^{i},\dot v^{i}\rp=\lp q^{i},\dot q^{i},v^{i},\dot v^{i}\rp$.

In order to describe $\alpha_{Q}$ is also necessary to define a tangent pairing. Given two manifolds $M$ and $N$, and a pairing $\bra\cdot,\cdot\ket:M\times N\rightarrow\mathbb{R}$ between them, the tangent pairing $\bra\cdot,\cdot\ket^{T}:TM\times TN\rightarrow\mathbb{R}$ is determined by
$$
\left\langle\frac{{\rm d}}{{\rm d}t}\Big |_{t=0}\hspace{1mm}\gamma\lp t\rp,\frac{{\rm d}}{{\rm d}t}\Big |_{t=0}\hspace{1mm}\delta\lp t\rp\right\rangle^{T}=\frac{{\rm d}}{{\rm d}t}\Big |_{t=0}\hspace{1mm}\langle\gamma\lp t\rp,\delta\lp t\rp\rangle
$$
where $\gamma:\mathbb{R}\rightarrow M$ and $\delta:\mathbb{R}\rightarrow N$.

Finally, we can define  $\alpha_{Q}$ as $\langle\alpha_{Q}\lp z\rp,w\rangle=\langle z,\kappa_{Q}\lp w\rp\rangle^{T}$, where $z\in TT^{*}Q$ and $w\in TTQ$. In local coordinates $\lp q^{i},p_{i}\rp$ for $T^{*}Q$ and $\lp q^{i},p_{i},\dot q^{i},\dot p_{i}\rp$ for $TT^{*}Q$, we have
$$
\alpha_{Q}\lp q^{i},p_{i},\dot q^{i},\dot p_{i}\rp=\lp q^{i},\dot q^{i},\dot p_{i},p_{i}\rp.
$$

 The  isomorphism $\beta_{Q}:TT^{*}Q\rightarrow T^{*}T^{*}Q$ is just given by $\beta_Q=\flat_{\omega_Q}$, where $\flat_{\omega_Q}$ is the
isomorphism defined by $\omega_Q$, that is, $\flat_{\omega_Q}(v)=i_v \omega_Q$.

 The Lagrangian dynamics is described by the Lagrangian submanifold ${\rm d}{\rm L}(TM)$ of $T^*TM$ where ${\rm L}\colon TM \rightarrow \mathbb{R}$ is the Lagrangian function, while the Hamiltonian formalism is described by the Lagrangian submanifold ${\rm d}H(T^*M)$ of $T^*TM$ where $H\colon T^*M \rightarrow \mathbb{R}$ is the corresponding Hamiltonian energy. In this paper we will extend the description of Lagrangian and Hamiltonian dynamics in terms of Lagrangian submanifolds and Morse families when controls are involved. Similar studies have already been developed by means of Lagrangian submanifolds in the literature for calculus of variations~\cite{2012FerLeMa}.

\subsection{Integrability algorithm}\label{SubSec:IntroIntegrAlgoritm}

Let $D$ be an implicit differential equation on a manifold $P$, that  is, a submanifold $D$ of $TP$.
In such a case, it is possible to construct an algorithm to extract the integrable part of $D$ in $P$ (see ~\cite{1995MMT}).

A curve $\gamma\colon I \subseteq \mathbb{R}\rightarrow P$
is called a \textit{solution of a differential equation} $D$ if $\dot{\gamma}(I)\subset D$.
An implicit differential equation $D$ is said to
be \textit{integrable at $v\in D$} if there
is a solution $\gamma\colon I \subseteq \mathbb{R}\rightarrow P$ such that
$\dot{\gamma}(0)=v$. An implicit differential equation $D$ is said to be integrable if it is integrable at each point
$v\in D$.

\begin{prop}\cite[Proposition 5]{1995MMT} Let $\tau_P\colon TP \rightarrow P$ be the canonical tangent bundle projection.
If $C=\tau_P(D)$ is a submanifold of $P$ and if the mapping
\begin{eqnarray*}
 \tau_P\colon & D& \rightarrow C\\
&v & \mapsto \tau_P(v)
\end{eqnarray*}
is a surjective submersion, then the condition $D\subset TC$ is sufficient for integrability
of the implicit differential equation $D$.
\end{prop}

In order to obtain the integrable part of an implicit differential equation we construct the following sequence of objects
\begin{equation*}
 (D^0,C^0,\tau^0),\; (D^1,C^1,\tau^1), \; \dots \; , \; (D^k,C^k,\tau^k), \; \dots
\end{equation*}
where
\begin{align}
 D^0&=D, \quad & C^0&= C=\tau_P(D), \quad &\tau^0&=\tau_P,\label{eq:StartIntegrAlg}\\
 D^k&=D^{k-1}\cap TC^{k-1}, \quad & C^k&= \tau_P(D^k), \quad &\tau^k&=\tau_P|_{D^k}\, .
\label{eq:kStepIntegrAlg}
\end{align}
For each $k$, it is assumed that the sets $C^k$ are submanifolds and that the mappings $\tau^k$
are surjective submersions. As the dimension of $P$ is finite, the sequence of implicit differential
equations $D^0$, $D^1,\dots , D^k,\dots$ stabilizes at some index $k=s$, that is,
$D^s=D^{s+1}$. The integrable implicit differential equation $D^s\subset TP$ is the integrable
part of $D$ and it could be empty.

This integrability algorithm, properly adapted, will be useful in Section~\ref{Sec:Integrab} to obtain the  absolutely continuous solutions
to optimal control problems which are the most common ones and correspond with piecewise constant controls.

\section{Introduction to optimal control problems} \label{Sec:IntroOCP}

Generally speaking, an optimal control problem from the differential geometric viewpoint is given by a vector field depending on parameters called controls, some boundary conditions and a cost function whose integral must be minimized or maximized.

\begin{defn} An \textit{optimal control problem} $({\mathcal C}, Q,\Gamma, {\rm L})$ is given by a control bundle
$\tau_{{\mathcal C},Q}\colon {\mathcal C} \rightarrow Q$, a vector field $\Gamma$ defined along the control bundle
projection $\tau_{{\mathcal C},Q}$, a cost function ${\rm L} \colon {\mathcal C} \rightarrow \mathbb{R}$ whose functional
must be minimized and some endpoint conditions or boundary conditions that must be satisfied at initial and/or final time.
\end{defn}

Remember that  such a vector field $\Gamma$ along $\tau_{{\mathcal C},Q}$ verifies $\tau_Q\circ \Gamma=\tau_{{\mathcal C},Q}$.
\begin{equation*}\xymatrix{&& TQ \ar[d]^{\txt{\small{$\tau_Q$}}} \\ {\mathcal C}
\ar[rr]^{\txt{\small{$\tau_{{\mathcal C},Q}$}}} \ar[rru]^{\txt{\small{$\Gamma$}}}& & Q}
\end{equation*}
Locally, $\dot{q}=\Gamma(q,u)$ and $\left(\tau_Q \circ \Gamma\right)(q,u)=\tau_{{\mathcal C},Q}(q,u)=q$.

From the optimal control data  $({\mathcal C}, Q, \Gamma, {\rm L})$ we construct the Pontryagin's hamiltonian $H:T^*Q\times_Q  {\mathcal C}\longrightarrow \R$ given by
\begin{equation}\label{pont}
H(\alpha_q, u_q)=\langle \alpha_q, \Gamma(u_q)\rangle-L(u_q)
\end{equation}
where $u_q\in C_q$ and $\alpha_q\in T_q^*Q$. 
In coordinates, $H(q^i, p_i, u^a)=p_j\Gamma^j(q^i, u^a)-L(q^i, u^a)$. This Hamiltonian is only valid for the so-called normal optimal solutions.
The usual technique to solve an optimal control problem is Pontryagin's Maximum Principle~\cite{Pontryagin}.

\begin{thm} \textbf{(Pontryagin's Maximum Principle, PMP)}
Let $U\subseteq \mathbb{R}^k$ and $(\gamma^*,u^*)\colon I \rightarrow {\mathcal C}=Q\times U$ be a normal solution of the
optimal control problem $({\mathcal C},Q,\Gamma,L)$. Then there exists $(\sigma^*, u^*)\colon
I \rightarrow T^*Q \times_Q  {\mathcal C}$ such that:
\begin{enumerate}
\item it is an integral curve of the Hamiltonian vector field $\Gamma_H$ with
given initial conditions in the states;
\item $\gamma^*=\pi_{Q} \circ \sigma^*$, with fiber
$\alpha^*(t)\in
T^*_{\gamma^*(t)} Q$ where $\pi_Q\colon T^*Q\rightarrow Q$ is the canonical projection to the cotangent bundle;
\item $H(\sigma^*(t),u^*(t))= \max_{u \in \overline{U}} H(\sigma^*(t),
u)$ almost everywhere, where $\overline{U}$ is the closure of $U$;
\item $\max_{u \in \overline{U}} H(\sigma^*(t),
u)$ is constant everywhere.
\end{enumerate} \label{thm:PMP}
\end{thm}

This theorem can be stated without the restriction of having a normal optimal solution. The abnormal solutions are associated with the Hamiltonian $H(\alpha_q, u_q)=\langle \alpha_q, \Gamma(u_q)\rangle$ that does not depend on the cost function.
See~\cite{Agrachev,2009BarberoMunoz,2007BressanPiccoli,Pontryagin} for more details on Pontryagin's Maximum Principle and optimal control theory.

This Principle provides us with a set of necessary conditions for optimality. In some optimal control problems the time interval is another unknown. If so, then the Hamiltonian function is zero along the optimal curve almost everywhere.

The condition of the maximization over the controls along the optimal solution is usually replaced by a more operative condition on the interior of $U$, that is,
\begin{equation}
\dfrac{\partial H}{\partial u^a}=0. \label{eq:HPartialU0}
\end{equation}

\section{Morse families and optimal control problems}\label{Sec:Morse}

In this section Morse families defined in Section~\ref{SubSec:LagranMorseIntro} are considered to describe geometrically the optimal control problems from the Hamiltonian viewpoint.

We first need to identify which is the candidate to be a Morse family in optimal control problems.  Without doubt the most natural candidate is the Pontryagin's Hamiltonian function $H\colon T^*Q\times_Q  {\mathcal C} \rightarrow \mathbb{R}$ defined in\ \eqref{pont}. According to Definition~\ref{Def:MorseFamily}, the submersion $\pi$ corresponds with $\pi\colon  T^*Q\times_Q  {\mathcal C} \rightarrow T^*Q$.
In adapted local coordinates, we have the following expressions that will be useful in the sequel:
\begin{eqnarray}
\ker {\rm T} \pi&=& \ds{\left\langle\frac{\partial}{\partial u^a}\right\rangle} , \nonumber \\
(T^*Q)_\pi &=& (\ker {\rm T}\pi)^0=\ds{\langle {\rm d}q^i, {\rm d}p_i \rangle},\nonumber  \\
{\rm d}H&=& \ds{\frac{\partial H}{\partial q^i} {\rm d}q^i+\frac{\partial H}{\partial p_i} {\rm d}p_i+\frac{\partial H}{\partial u^a} {\rm d}u^a},\nonumber  \\
(T^*Q)_\pi \cap {\rm d}H(T^*Q\times_Q  {\mathcal C}) &=& \ds{\left\{ (q^i,p_i,u^a,A_i,B^i,C_a)\in T^*(T^*Q\times_Q  {\mathcal C}) \, | \, A_i=\frac{\partial H}{\partial q^i}, \right. } \nonumber \\&& \quad \quad \ds{\left.\;  B^i =\frac{\partial H}{\partial p_i}, \; C_a=\frac{\partial H}{\partial u^a}=0 \right\}}. \label{eq:MorseIntersection}
\end{eqnarray}
From Definition~\ref{Def:MorseFamily} we know that the  Pontryagin's Hamiltonian   $H\colon T^*Q\times_Q  {\mathcal C}\rightarrow \mathbb{R}$ defines a Morse family  if and only if the submanifolds $ (T^*Q)_\pi$ and ${\rm d}H(T^*Q\times_Q  {\mathcal C})$ of $T^*(T^*Q\times_Q  {\mathcal C})$ are transverse.
%

From Section~\ref{SubSec:LagranMorseIntro} we obtain the following result.
\begin{prop} Pontryagin's Hamiltonian   $H\colon T^*Q\times_Q  {\mathcal C}\rightarrow \mathbb{R}$ defines a Morse family if and only if  the matrix
\begin{equation}
{\rm D}_{(q,p,u)} \begin{pmatrix} \dfrac{\partial H}{\partial u}\end{pmatrix}=\begin{pmatrix}\dfrac{\partial^2 H}{\partial q^i \partial u^a} & \dfrac{\partial^2 H}{\partial p_i \partial u^a} & \dfrac{\partial^2 H}{\partial u^a \partial u^b}
\end{pmatrix}_{(q,p,u)}  \label{Eq:SecondDerivativeHu}
\end{equation}
has maximum rank for all $(q,p,u)\in T^*Q\times_Q {\mathcal C}$.
\begin{proof}
Let us check when the conormal bundle $(T^* Q)_\pi$ and $ {\rm d}H(T^*Q\times_Q  {\mathcal C})$ are transverse in $T^*(T^*Q\times_Q {\mathcal C})$, that is,  when condition~\eqref{Eq:Transverse},
\begin{equation*} {\rm T}_\alpha( {\rm d}H(T^*Q\times_Q  {\mathcal C}))+T_\alpha (T^*Q)_\pi=T_\alpha (T^*(T^*Q\times_Q  {\mathcal C})),
\end{equation*}
is satisfied for all $\alpha\in (T^* Q)_\pi\cap {\rm d}H(T^*Q\times_Q  {\mathcal C})$. The local generators of the subspaces on the left-hand side of the above equation, given in local
coordinates $(q,p,u,A,B,C)$, are the following ones:
\begin{eqnarray}
T_\alpha (T^*Q)_\pi&=&{\rm span}_{\mathbb{R}} \left\{  \dfrac{\partial}{\partial q^i},\,
\dfrac{\partial}{\partial p_i},\, \dfrac{\partial}{\partial u^a},\, \dfrac{\partial}{\partial A_i},\, \dfrac{\partial}{\partial B_i}  \right\}_\alpha, \label{eq:Prop15-1}\\
T_\alpha {\rm d}H(T^*Q\times_Q  {\mathcal C} )&=& {\rm span}_{\mathbb{R}} \left\{ \dfrac{\partial}{\partial q^i}+\dfrac{\partial^2 H}{\partial q^i \partial q^j}\, \dfrac{\partial}{\partial A_j}+\dfrac{\partial^2 H}{\partial q^i \partial p_j}\, \dfrac{\partial}{\partial B_j}+\dfrac{\partial^2 H}{\partial q^i \partial u^a}\, \dfrac{\partial}{\partial C_a}, \right.
\nonumber \\
&& \qquad \qquad  \dfrac{\partial}{\partial p_i}+\dfrac{\partial^2 H}{\partial p_i \partial q^j}\, \dfrac{\partial}{\partial A_j}+\dfrac{\partial^2 H}{\partial p_i \partial u^a}\, \dfrac{\partial}{\partial C_a},
\label{eq:Prop15-2} \\ & & \left.   \qquad \qquad \dfrac{\partial}{\partial u^a}+\dfrac{\partial^2 H}{\partial u^a \partial q^j}\, \dfrac{\partial}{\partial A_j}+\dfrac{\partial^2 H}{\partial u^a \partial p_j}\, \dfrac{\partial}{\partial B_j}+\dfrac{\partial^2 H}{\partial u^a \partial u^b}\, \dfrac{\partial}{\partial C_b}
\right\}_\alpha. \nonumber
\end{eqnarray}
Note that these two tangent spaces satisfy the condition~\eqref{Eq:Transverse} if the matrix
\begin{equation*}
\begin{pmatrix}
\dfrac{\partial^2 H}{\partial q^i \partial u^a} & \dfrac{\partial^2 H}{\partial p_i \partial u^a} & \dfrac{\partial^2 H}{\partial u^a \partial u^b}
\end{pmatrix}
\end{equation*}
has maximum rank. Only if this condition is satisfied the Hamiltonian function is a Morse family. Otherwise, there is only a strict inclusion ${\rm T}_\alpha( {\rm d}H(T^*Q\times_Q  {\mathcal C}))+T_\alpha (T^*Q)_\pi\subsetneq T_\alpha (T^*(T^*Q\times_Q  {\mathcal C}))$.
\end{proof} \label{Prop:PHMorse}
\end{prop}
\begin{remark}
In this paper we restrict ourselves to examples where the optimal control problem is defined by a Morse family given by the Pontryagin's Hamiltonian
$H: T^*Q\times_Q {\mathcal C}\rightarrow \R$. However there are many interesting situations where this condition is not fulfilled. For instance, the optimal control problem given by the control equation $\dot{q}=u^2$ and cost function $L\equiv 1$ does not satisfy the Morse family condition
(\ref{Eq:SecondDerivativeHu}) for all points $(q,p,u)$ in $T^*Q\times_Q {\mathcal C}$ since $H(q, p, u)=pu^2-1$ and
\[
 D_{(q,p,u)}\left(\frac{\partial H}{\partial u}\right)=\left( \begin{array}{lll} 0& 2u& 2p  \end{array}\right).
 \]
Note that $H$ does not define a Morse family at the points $(q, 0, 0)$ because at these points the matrix $D_{(q,p,u)}\left(\frac{\partial H}{\partial u}\right)$ does not have maximum rank. In fact, it is also possible to check in this example
that the intersection of ${\rm d}H(T^*Q\times_Q {\mathcal C})$ and $(T^*Q)_{\pi}$ is not clean at all the points in $T^*Q\times_Q {\mathcal C}$.
\end{remark}

If $H: T^*Q\times_Q {\mathcal C}\rightarrow \R$ defines a Morse family, then $\Sigma_H={\rm d}H(T^*Q\times_Q  {\mathcal C}) \cap  (T^*Q)_\pi$ is a submanifold of $T^*(T^*Q\times_Q  {\mathcal C})$.
The Pontryagin's Hamiltonian is called the generating function  for the Lagrangian submanifold
 ${\mathcal L}_H={\rm pr}(\Sigma_H)$  of
$(T^*T^*Q, \omega_{T^*Q})$, where ${\rm pr}\colon T^*(T^*Q
\times_Q {\mathcal C}) \rightarrow T^*T^*Q$ is the natural projection onto the first factor, see Proposition
\ref{prop:img:j:lagrangian}. Locally,
\begin{align}
{\mathcal L}_H=\left\{(q, p, P_q, P_p)\in T^*(T^*Q)\; |\; \exists \; u \hbox{ such that }\right.  & (q,u)\in {\mathcal C}, \;  P_q=\dfrac{\partial H}{\partial q}(q,p,u), P_p=\dfrac{\partial H}{\partial p}(q,p,u)  \nonumber 
\\
& \left.   \hbox{ and }\frac{\partial H}{\partial u}(q,p,u)=0\right\} \label{eq:DefLH}
\end{align}

Depending on which submatrix of ${\rm D}_{(q,p,u)} \begin{pmatrix} \dfrac{\partial H}{\partial u}\end{pmatrix}$
gives the maximum rank in condition~\eqref{Eq:SecondDerivativeHu}, there are different kinds of optimal control problems associated with different relative positions of ${\mathcal L}_H$ with respect to the canonical projection $\pi_{T^*Q}: T^*T^*Q\rightarrow T^*Q$.

\begin{defn}\label{def:regular:opc}
We say that the optimal control problem is \textit{regular} if the restriction $(\pi_{T^*Q})_{|{\mathcal L}_H}: {\mathcal L}_H\rightarrow T^*Q$ is a local diffeomorphism. \label{Defn:RegOCP}
\end{defn}
If the optimal control problem is regular, the submanifold ${\mathcal L}_H$ is transverse to the fibers of $\pi_{T^*Q}: T^*T^*Q\rightarrow T^*Q$.
Locally, the regularity condition is satisfied  when the
matrix ${\rm D}_{(q,p,u)} \begin{pmatrix} \dfrac{\partial H}{\partial u}\end{pmatrix}$ has maximum rank because
$\dfrac{\partial^2H}{\partial u^a \partial u^b}$ has maximum rank. If so,  for each $(q,p)\in T^*Q$ there exists a function $u^*=u^*(q,p)$ by the Implicit Function Theorem defined in a neighbourhood of $(q,p)$ such that
\begin{equation*}
\left(q,p,\dfrac{\partial H}{\partial q^i}(q,p,u^*(q,p)), \dfrac{\partial H}{\partial p_i}(q,p,u^*(q,p)),0  \right)\in {\mathcal L}_H.
\end{equation*}
In other words, ${\rm d}H_{u^*}(T^*Q)={\mathcal L}_H$
where $H_{u^*}(q,p)=H(q, p, u^*(q, p))$.

If the condition in Definition~\ref{Defn:RegOCP} is not fulfilled, then the optimal control problem is \textit{singular}. This singularity may appear because of many different reasons such as $\pi_{T^*Q}( {\mathcal L}_{H})\varsubsetneq T^*Q$ or/and
there are points $\mu\in {\mathcal L}_H$ where
\[
{\rm T}_{\mu}\pi_Q: T_\mu {\mathcal L}_{H}\longrightarrow T_{\alpha} T^*Q, \qquad \alpha=\pi_{T^*Q}(\mu)
\]
is not surjective, that is, $\mu$ is a \textit{caustic point}~\cite{1990Arnold,2012ArnoldGV}.

Let us consider now simple examples to show different situations.
\begin{example}
Consider the optimal control problem given by the control equation $\dot{q}=u$ and $L=\frac{1}{6}u^3$ with $u\in \R$. Pontryagin's Hamiltonian $H(q,p,u)=pu-\dfrac{u^3}{6}$ is obviously a Morse family because the matrix
\[
 D_{(q,p,u)}\left(\frac{\partial H}{\partial u}\right)=\left( \begin{array}{lll} 0& 1& -u  \end{array}\right)
 \]
has maximum rank. But the optimal control problem is not regular if $u=0$. 
The submanifold ${\mathcal L}_H$ in~\eqref{eq:DefLH} is given by 
\begin{eqnarray*}
{\mathcal L}_H&=&\left\{(q, p, P_q, P_p)\in T^*(T^*\R)\equiv \R^4\; |\; P_q=0, P_p=u, p=\frac{1}{2}u^2\right\}\\
              &=&\left\{(x, \frac{1}{2}\lambda^2, 0, \lambda)\, | \lambda\in \R\right\}
\end{eqnarray*}
Observe that the points $(x, 0, 0, 0)$ are caustic and moreover the
$\pi_{T^*Q}( {\mathcal L}_{H})=\mathbb{R}\times [0, \infty)\varsubsetneq \R^2$.  From the local expression of ${\mathcal L}_H$ it is also clear that $(\pi_{T^*Q})_{|{\mathcal L}_H}\colon {\mathcal L}_H \rightarrow T^*Q$ 
is a local diffeomorphism as long as $u\neq 0$. 
\end{example}

\begin{example}[Singular optimal control problem for LQ systems]
Let us consider a general example of singular optimal control problems, see~\cite{2002DelgadoIbort, 2009DelgadoIbort,1999BelgasCDC} for more details.
These problems are characterized by the impossibility to determine the controls with first-order necessary conditions for optimality. 

Consider the control system
\begin{equation*}
\dot{q}(t)=Aq(t)+Bu(t)=\Gamma(q,u)
\end{equation*}
with cost function
\begin{equation*} {\rm L}(q,u)=\dfrac{1}{2}u^{\rm T}P u +q^{\rm T}Q u +\dfrac{1}{2}q^{\rm T}R q,
\end{equation*}
where all matrices $A$, $B$, $P$, $Q$ and $R$ are constant and of suitable size.

Pontryagin's Hamiltonian is given by $H(q,p,u)=p^{\rm T} \left(Aq+Bu\right)-{\rm L}(q,u)$.
As described above, the candidate to be a Morse family is the Pontryagin's Hamiltonian.
Let $\pi\colon T^*Q\times_Q {\mathcal C}\rightarrow T^* Q$, the local expressions for the elements in equation~\eqref{eq:MorseIntersection} are:
\begin{eqnarray*}
\ker {\rm T} \pi&=& \left\{ (q,p,u,V_q,V_p,V_u)\in T( T^*Q\times_Q {\mathcal C}) \; | \; V_q=0, \; V_p=0\right\},\\
(T^*Q)_\pi &=& (\ker {\rm T}\pi)^0=\left\{ (q,p,u,A,B,C)\in T^*(T^*Q\times_Q  {\mathcal C}) \, | \,  C=0\right\},\\
{\rm d}H&=& \ds{\left(p^{\rm T}A-Q u-Rq\right)_i {\rm d}q^i+\left( Aq+Bu\right)^i {\rm d}p_i+\left( p^{\rm T}B-Pu-q^{\rm T}Q\right)_a {\rm d}u^a},\\
(T^*Q)_\pi \cap {\rm d}H(T^*Q\times_Q  {\mathcal C}) &=& \ds{\left\{ (q,p,u,A,B,C)\in T^*(T^*Q\times_Q  {\mathcal C}) \, | \, A_i=\left(p^{\rm T}A-Q u-Rq\right)_i, \right. } \\&& \quad \quad \ds{\left.\;  B^i =\left( Aq+Bu\right)^i, \; C_a=\left( p^{\rm T}B-Pu-q^{\rm T}Q\right)_a=0 \right\}}.
\end{eqnarray*}
Let us study the rank of the matrix in condition~\eqref{Eq:SecondDerivativeHu},
\begin{equation} {\rm D}_{(q,p,u)} \left( \dfrac{\partial H}{\partial u^a} \right)=
\begin{pmatrix}
\dfrac{\partial^2 H}{\partial q^i \partial u^a} & \dfrac{\partial^2 H}{\partial p_i \partial u^a} & \dfrac{\partial^2 H}{\partial u^a \partial u^b}
\end{pmatrix}=\begin{pmatrix} -Q & B & -P \end{pmatrix}, \label{eq:ExampleMorseHessianControl}
\end{equation}
to decide if the Pontryagin's Hamiltonian is a Morse family and if it defines a regular or singular optimal control problem.
If the above matrix has maximum rank, Pontryagin's Hamiltonian will be a Morse family associated with either a regular or singular optimal control problem depending on what submatrices of~\eqref{eq:ExampleMorseHessianControl} have maximum rank.
This optimal control problem is regular if $P$ has maximum rank.
Otherwise, the optimal control problem is singular.
\end{example}

\section{The integrability algorithm for optimal control problems}\label{Sec:Integrab}
Given an optimal control problem we have defined in the previous section the subset ${\mathcal L}_H$ in~\eqref{eq:DefLH} which is an immersed Lagrangian submanifold, as long as the Pontryagin's Hamiltonian defines a Morse family.
Independently of the manifold structure of ${\mathcal L}_H$, we always have an optimal control problem associated with the subset
\[
 {\rm pr}\left( (T^*Q)_\pi \cap {\rm d}H(T^*Q\times_Q  {\mathcal C}) \right)={\mathcal L}_H \subset T^*T^*Q.
 \]
Thus, using the musical  isomorphism
\[
\sharp_{T^*Q}: T^*T^*Q\longrightarrow TT^*Q
\]
defined by the canonical symplectic 2-form $\omega_Q$ on $T^*Q$, we obtain a subset 
\begin{equation}{\mathcal D}_H=\sharp_{T^*Q}({\mathcal L}_H)
\label{eq:DH}
\end{equation}
 of $TT^*Q$. In this space 
it is possible to apply the integrability algorithm described in Section~\ref{SubSec:IntroIntegrAlgoritm}. Here we will not apply the algorithm, but adapt it to 
take into account all the peculiarities of an optimal control problem.
For instance, in an optimal control problem the controls are always measurable and bounded, and in general piecewise constant. As a result the
curves on $T^*Q$ satisfying the necessary conditions for optimality in Theorem~\ref{thm:PMP} are absolutely continuous. Hence, they
are differentiable almost everywhere and satisfy most of the necessary conditions almost everywhere.

Unless otherwise stated, in this section the optimal control problems are assumed to be regular. Hence, for
each point in $T^*Q$, we can describe the control function $u$ in terms of the states and momenta from the equation $\dfrac{\partial H}{\partial u}(q,p,u)=0$.

A curve $\sigma$ on $T^*Q\times_Q {\mathcal C}$ satisfies the necessary conditions for optimality in Theorem~\ref{thm:PMP}
if the curve $\gamma\equiv  {\rm pr}_1\circ {\sigma}\colon I\subset \mathbb{R} \rightarrow T^*Q$, where ${\rm pr}_1\colon T^*Q\times_Q {\mathcal C} \rightarrow T^*Q$, satisfies that  $\dot{\gamma}(I)\subset \sharp_{T^*Q}({\mathcal L}_H)$ almost everywhere.

 Note that ${\mathcal D}_H$
is a subset of  $TT^*Q$ with the following local expression
\begin{align*}
{\mathcal D}_H=\left\{(q, p, V_q, V_p)\in T(T^*Q)\; |\; \exists \;  u \hbox{ such that }  (q,u)\in {\mathcal C}, \; \right. & V_q=\frac{\partial H}{\partial p}(q,p,u), V_p=-\frac{\partial H}{\partial q}(q,p,u) \\
& \left.\hbox{ and }\frac{\partial H}{\partial u}(q,p,u)=0\right\}\, .
\end{align*}

%
%

In order to include the possibility that the solutions are piecewise $C^1$ curves we need to modify the integrability algorithm described in Section~\ref{SubSec:IntroIntegrAlgoritm}.

Let $\tau_{T^*Q}\colon TT^*Q\rightarrow T^*Q$ be
the canonical tangent bundle projection and consider the subset ${\mathcal D}_H\times_{\tau_{T^*Q}({\mathcal D}_H)}{\mathcal D}_H$ of the Whitney sum $T(T^*Q)\oplus T(T^*Q)$.

A piecewise $C^1$ curve $\gamma\colon I \subseteq \mathbb{R}\rightarrow T^*Q$
is called a \textit{solution of the optimal control problem } ${\mathcal D}_H\times_{\tau_{T^*Q}({\mathcal D}_H)}{\mathcal D}_H$ if $(\dot{\gamma}^-(t),\dot{\gamma}^+(t))\in {\mathcal D}_H\times_{\tau_{T^*Q}({\mathcal D}_H)}{\mathcal D}_H$ for all $t\in I$.
The implicit differential equation determined by ${\mathcal D}_H\times_{\tau_{T^*Q}({\mathcal D}_H)}{\mathcal D}_H$ is said to
be \textit{integrable at $(v_1,v_2)\in {\mathcal D}_H\times_{\tau_{T^*Q}({\mathcal D}_H)}{\mathcal D}_H$} if there
is a solution $\gamma\colon I \subseteq \mathbb{R}\rightarrow T^*Q$ such that
$(\dot{\gamma}^-(0),\dot{\gamma}^+(0))=(v_1,v_2)$.
Note that this notion of integrability includes the points of both
continuity and discontinuity
of the derivative.

Now we can adapt the integrability algorithm in Section~\ref{SubSec:IntroIntegrAlgoritm} by defining
the following sequence of objects
\begin{equation*}
 D^0={\mathcal D}_H\times_{\tau_{T^*Q}({\mathcal D}_H)}{\mathcal D}_H, \quad C^0 =\tau_{T^*Q}({\mathcal D}_H\times_{\tau_{T^*Q}({\mathcal D}_H)}{\mathcal D}_H),\quad
\tau^0=\tau_{T^*Q}|_{{\mathcal D}_H\times_{\tau_{T^*Q}({\mathcal D}_H)}{\mathcal D}_H},
\end{equation*}
and
\begin{eqnarray*}
D^k&=&\left\{(v_1,v_2)\in D^{k-1} \; | \; \exists \, \gamma\colon I \rightarrow C^{k-1}, \quad
\dot{\gamma}^{-}(0)=v_1, \quad    \dot{\gamma}^{+}(0)=v_2\right\}\,,\\
C^k&=&\tau_{T^*Q}(D^k)\, ,\\
\tau^k&=& \tau_{T^*Q}|_{D^k}\colon D^k\rightarrow T^*Q.
\end{eqnarray*}
At some step $s$ the algorithm stabilizes, that is, $D^s=D^{s+1}$. The integrable
implicit differential equation $D^s\subset T(T^*Q)\oplus T(T^*Q)$ is the integrable
part, maybe empty, of ${\mathcal D}_H\times_{\tau_{T^*Q}({\mathcal D}_H)}{\mathcal D}_H$.

\example  Consider the following optimal control problem
\begin{eqnarray*}
{\rm min} \, J&=&\int^{t_f}_{t_0}(u(t)^2-1)^2{\rm d}t,\\
\mbox{subject to} & \dot{q}=u, &  U\subseteq \mathbb{R}, \\
& t_0=0, & q_0=1,\\
& t_f=1, & q_f=1.
\end{eqnarray*}
 Pontryagin's Hamiltonian associated with this problem is $H(q,p,u)=pu-(u^2-1)^2$. Obviously, this Hamiltonian defines a Morse family because
the matrix
\begin{equation*}
D_{(q,p,u)}\left( \dfrac{\partial H}{\partial u} \right)=\left( 0 \; 1 \; -12u^2+4 \right)
\end{equation*}
has maximum rank for all $(q,p,u)$ in $T^*\mathbb{R}\times_{\mathbb{R}}{\mathcal C}$, see Proposition~\ref{Prop:PHMorse}.

The set ${\mathcal D}_H$ of $TT^*Q$ in this example has the following local expression:
\begin{eqnarray*}
{\mathcal D}_H &=&\left\{ (q,p,V_q,V_p)\in T(T^* Q)\; | \; \exists \; u \mbox{ such that }(q,u)\in {\mathcal C}, \;  V_q=\dfrac{\partial H}{\partial p}(q,p,u), \;  V_p=-\dfrac{\partial H}{\partial q}(q,p,u),\; \right.
\\ && \qquad \left. \dfrac{\partial H}{\partial u}(q,p,u)=0  \right\}\\
&=&  \left\{ (q,p,V_q,V_p)\in T(T^* Q)\; | \; \exists \; u \mbox{ such that } (q,u)\in {\mathcal C}, \;  V_q=u, \;  V_p=0,\; p-4u(u^2-1)=0  \right\}\\
&=&  \left\{ (x,4u(u^2-1),u,0)\in T(T^* Q)\; | \; \exists \; u \mbox{ such that $u$ is constant}   \right\}.\\
\end{eqnarray*}
The last equality comes from the fact that $\dot{p}=0$. Hence $p$ is a constant momenta equal to $4u(u^2-1)$. Thus the control $u$ must also be constant.

The starting point of the integrability algorithm is given by
\begin{eqnarray*}
D^0&=&{\mathcal D}_H\times_{\tau_{T^*Q}({\mathcal D}_H)} {\mathcal D}_H\\&=&\left\{ (x,4u_1(u^2_1-1),u_1,0)\in T(T^* Q) \; |\; u_1\in U\right\}\times_{\tau_{T^*Q}({\mathcal D}_H)}\\&& \quad
\left\{ (x,4u_2(u^2_2-1),u_2,0)\in T(T^* Q) \; |\; u_2\in U \right\},\\
C^0&=&\tau_{T^*Q}({\mathcal D}_H\times_{\tau_{T^*Q}({\mathcal D}_H)} {\mathcal D}_H)=\left\{(x,4u_1(u^2_1-1))\; | \; u_1\in U\right\},\\
\tau^0&=&\tau_{T^*Q}|_{{\mathcal D}_H\times_{\tau_{T^*Q}({\mathcal D}_H)}{\mathcal D}_H}.
\end{eqnarray*}
As a consequence of being a fibered product over $T^* Q$, $4u_1(u^2_1-1)=4u_2(u^2_2-1)$ so that $C^0$ is well defined.
That  condition is fulfilled when $u_1=u_2$ or if there is a switching time in the trajectory where the control changes its value, that is, $u_1\neq u_2$.
However, if $u_1\neq u_2$, the equality $4u_1(u^2_1-1)=4u_2(u^2_2-1)$ is satisfied if and only if $u_1,u_2\in \{-1,0,1\}$. Since the trajectories are straight lines, the controls can only change from $1$ to $-1$ or vice versa in order to verify the endpoint conditions if only one switching time is allowed.
For the case $u_1=u_2$, the endpoint conditions are only satisfied if $u_1=u_2=0$.
Under the assumption of not having more than one change in the control, the next step in the integrability algorithm is the following one:
\begin{eqnarray*}
D^1&=&\left\{ ((x,0,u_1,0),(x,0,u_2,0))\in T(T^* Q)\times _{\tau_{T^*Q}({\mathcal D}_H)}T(T^*Q) \; | \; \{u_1=\pm 1, u_2=\mp 1\} \right. \\
&& \qquad \left. \mbox{ or } \{u_1=u_2=0\} \right\},\\
C^1&=&\left\{(x,0)\right\},\\
\tau^1&=&\tau_{T^*Q}|_{D^1}\colon D^1 \rightarrow T^*Q.
\end{eqnarray*}
Here the algorithm stabilizes and the solutions are
\begin{itemize}\item $
\gamma_1(t)=\left\{ \begin{array}{rcl} u_1t+1 && 0\leq t \leq t_1, \\
u_2(t-t_f)+1 && t_1\leq t \leq t_f, \end{array}\right.
$
where $u_1=\pm 1$, $u_2=\mp 1$ such that $\gamma_1^-(t_1)=\gamma_1^+(t_1)$, that is, $u_1t_1+1=u_2(t_1-t_f)+1$. Equivalently, $t_1=t_f/2$.
\item $\gamma_2(t)=1$ for $u_1=u_2=0$.
\end{itemize}
If we compute the value of the functional to be minimized along these two families of trajectories, we obtain that
\begin{eqnarray*}
 {\rm J}(\gamma_1)&=&0,\\
{\rm J}(\gamma_2)&=&1.
\end{eqnarray*}
Thus the optimal solution comes from curves like $\gamma_1$.

\section{Towards a Lagrangian formulation of an optimal control problem}\label{Sec:LagranOCP}

After introducing in Section~\ref{Sec:Morse} the  Hamiltonian
approach to optimal control problems by Morse families, we are going to consider a Lagrangian approach to deal with optimal control problems by means of 
Lagrangian submanifolds of $T^*(TQ)$ instead of $T^*(T^*Q)$. In order to do this we make use of the Tulczyjew diffeomorphism which defines a diffeomorphism between $T^*TQ$ and $T^*T^*Q$, see Section \ref{Sec:Tulczy}.

Given $({\mathcal C}, Q, \Gamma, {\rm L})$ an optimal control problem,  construct the following subset of $T^*TQ$:
\begin{equation}
{\mathcal L}_{{\mathcal C},{\rm L}}=\{\mu \in T^*TQ \; | \; \Gamma^* \mu={\rm d}{\rm L} \}.
\label{eq:GammaCL}
\end{equation}
In optimal control theory there are two different kinds of solutions: the normal ones, that will be obtained from ${\mathcal L}_{{\mathcal C},{\rm L}}$, and the abnormal ones. The abnormality is characterized by the no dependence on the cost function at a first stage. These abnormal curves can be included in our description by considering the following subset in $T^*TQ$:
\begin{equation*}
{\mathcal L}_{{\mathcal C},abn}=\{\mu \in T^*TQ \; | \; \Gamma^* \mu=0 \}.
\label{eq:GammaCAbn}
\end{equation*}
From now on, we just consider the characterization of normal solutions, but analogous constructions can be considered for abnormal solutions.

In local coordinates $(q^i,\dot{q}^i,u^a)$ for $TQ\times_Q {\mathcal C}$, the elements in~\eqref{eq:GammaCL} are written as follows:
\begin{equation*} \begin{array}{rcl}
\mu&=&a_i{\rm d}q^i+b_i {\rm d} \dot{q}^i,\\
\Gamma^*\mu&=&a_i {\rm d}q^i+b_i\ds{\frac{\partial \Gamma^i}{\partial q^j} \, {\rm d} q^j+b_i\frac{\partial \Gamma^i}{\partial u^a}\, {\rm d}u^a}, \\
{\rm d}{\rm L}&=& \ds{\frac{\partial L}{\partial q^i} \, {\rm d}q^i+ \frac{\partial L}{\partial u^a} \, {\rm d}u^a}.
\end{array}
\end{equation*}
Hence the elements in ${\mathcal L}_{{\mathcal C},{\rm L}}$ satisfy the following equations:
\begin{eqnarray}
a_i+b_j\ds{\frac{\partial \Gamma^j}{\partial q^i}}&=&\ds{\frac{\partial {\rm L}}{\partial q^i}},\label{eq:LocalGammaCL-1}\\
b_i \ds{\frac{\partial \Gamma^i}{\partial u^a}}&=&\ds{\frac{\partial {\rm L}}{\partial u^a}},\label{eq:LocalGammaCL-2}
\end{eqnarray}
restricted to $\dot{q}^i=\Gamma^i(q,u)$. 
Condition (\ref{eq:LocalGammaCL-2}) may introduce additional restrictions in the control space $C$ as shown in the following lemma.

\begin{lemma}
Define the subset $\widetilde{{\mathcal C}}\subseteq {\mathcal C}$ as
\[
\widetilde{{\mathcal C}}=\{ c\in {\mathcal C}\; |\; {\rm d}{\rm L}(V_c)=0, \; \forall \; V_c\in \ker {\rm T}_c\Gamma\},
\]
then $\pi_{TQ}({\mathcal L}_{{\mathcal C},{\rm L}})\subseteq \Gamma (\widetilde{{\mathcal C}})$.
\proof
Let $v_q\in \pi_{TQ}({\mathcal L}_{{\mathcal C},{\rm L}})$, there exists $\mu_{v_q}$ in ${\mathcal L}_{{\mathcal C},{\rm L}}$ such that
\[
\Gamma^*(\mu)={\rm d}{\rm L} \mbox{ at } v_q, \mbox{ that is, }\langle\mu_{v_q}, {\rm T}_c\Gamma(V_c)\rangle={\rm dL}(V_c),
\]
for all $V_c\in {\rm T}_c{\mathcal C}$ such that $\Gamma(c)=v_q$. 
In particular, for all  $V_c\in \ker {\rm T}_c\Gamma$ we have that
$dL(V_c)=0$. Hence, $c\in \widetilde{{\mathcal C}}$ and  $v_q\in \Gamma(\widetilde{{\mathcal C}})$. 
\qed \label{Lemma:TildeC}
\end{lemma}   
Observe that if $\Gamma$ is an immersion, then $\widetilde{{\mathcal C}}={\mathcal C}$. 

Assume in the sequel that $\widetilde{{\mathcal C}}$ is transverse to the fibers of $\Gamma$, that is
\[
{\rm T}\widetilde{{\mathcal C}}+ \ker {\rm T}\Gamma={\rm T} {\mathcal C}.
\]
Then it is easy to prove the following result.
\begin{prop}\label{propo-3} Let $\widetilde{{\mathcal C}}=\{ c\in {\mathcal C}\; |\; {\rm d}{\rm L}(V_c)=0, \; \forall \; V_c\in \ker {\rm T}_c\Gamma\},$ we have
\[
{\mathcal L}_{{\mathcal C},{\rm L}}={\mathcal L}_{\widetilde{\mathcal C},{\rm L}_{|\widetilde{\mathcal C}}}.
\]
\end{prop}
Observe that the restriction of $\tau_{{\mathcal C}, Q}\colon {\mathcal C} \rightarrow Q$ to  $\widetilde{\mathcal C}$ is not necessarily the entire space $Q$. However, Proposition~\ref{propo-3} 
guarantees we can work indistinctly with  ${\mathcal L}_{{\mathcal C},{\rm L}}$ and ${\mathcal L}_{\widetilde{\mathcal C},{\rm L}_{|\widetilde{\mathcal C}}}$.

\begin{example}
As a toy example consider the optimal control system given by $\dot{x}=u_1+u_2$ and cost function $L=x(u_1-u_2)$. In this case, 
\begin{eqnarray*}
{\mathcal L}_{{\mathcal C},{\rm L}}&=&\{(x,v, a, b)\in \R^4\; |\; \exists \;  (u_1, u_2) \in \mathbb{R}^2 \hbox{ such that }
v=u_1+u_2, a=u_1-u_2, b=x, b=-x\}\\
&=&\{(0, r_1, r_2, 0)\;|\; (r_1, r_2)\in \R^2\}
\end{eqnarray*}
On the other hand, $\widetilde{\mathcal C}=\{(x, u_1, u_2)\in {\mathcal C} | x=0,\;  u_1=u_2\}$. Note that $\pi_{TQ}({\mathcal L}_{{\mathcal C},{\rm L}})\subseteq \Gamma (\widetilde{{\mathcal C}})$ as proved
in Lemma~\ref{Lemma:TildeC}.
\end{example}

Let us consider the Pontryagin's Hamiltonian function $H: T^*Q\times_Q  {\mathcal C}\longrightarrow \R$ with local expression
\begin{equation*}
H(q,p,u)=p_i\Gamma^i(q,u)-{\rm L}(q,u).
\end{equation*}
Then the equations in~\eqref{eq:LocalGammaCL-1} and~\eqref{eq:LocalGammaCL-2} can be rewritten as follows:
\begin{equation}
\begin{array}{rcl}
a_i&=&-\ds{\frac{\partial H}{\partial q^i}(q,b,u)},\\
0&=&\ds{\frac{\partial H}{\partial u^a}(q,b,u).}
\end{array}
\end{equation}
restricted to $\dot{q}^i=\ds{\frac{\partial H}{\partial p_i}(q,b,u)}=\Gamma^i(q, u)$. Locally,
\begin{equation*}
{\mathcal L}_{{\mathcal C},L}=\left\{ (q^i,\dot{q}^i,a_i,b_i)\in T^*TQ  \; | \; \exists\; u \hbox{ such that } (q,u)\in {\mathcal C}, \, \dot{q}^i=\frac{\partial H}{\partial p_i}, \, a_i=-\frac{\partial H}{\partial q^i}, \, \frac{\partial H}{\partial u^a}=0 \right\}. \label{eq:LocalGammaCLset}
\end{equation*}
These are the first order necessary conditions for optimality in Pontryagin's Maximum Principle, providing that the controls
take value in an open set so that the maximization condition of the Hamiltonian function over the controls is replaced
by the weaker condition~\eqref{eq:HPartialU0}.

By means of the Tulczyjew diffeomorphism 
\begin{equation*}
\begin{array}{crcl}
\alpha_Q \colon & TT^*Q & \longrightarrow & T^*TQ\\
& (q^i,p_i,\dot{q}^i,\dot{p}_i) & \longmapsto & (q^i,\dot{q}^i,\dot{p}_i,p_i),
\end{array}
\label{eq:TulczyDiffeo}
\end{equation*}
and~\eqref{eq:DH} it is easy to show that the subsets ${\mathcal L}_H$ defined in~\eqref{eq:DefLH}
and ${\mathcal L}_{{\mathcal C},{\rm L}}$ satisfy that
\[
{\mathcal L}_{{\mathcal C},{\rm L}}=\alpha_Q({\mathcal D}_H)=(\alpha_Q\circ \sharp_{T^*Q})({\mathcal L}_H)\; .
\]
From Proposition~\ref{propo-3} we have that ${\mathcal L}_{\widetilde{\mathcal C},{\rm L}_{|\widetilde{\mathcal C}}}={\mathcal L}_{{\mathcal C},{\rm L}}=\alpha_Q({\mathcal D}_H)$.
The above equalities imply the following result.
\begin{prop} \label{Prop:LagrSubmanifold}  Let $({\mathcal C},Q,\Gamma, {\rm L})$ be an optimal control problem such that the associated Pontryagin's Hamiltonian is a Morse
 family, then  ${\mathcal L}_{{\mathcal C},{\rm L}}$ is an
immersed Lagrangian submanifold of the symplectic manifold $(T^*TQ,\omega_{TQ})$.
\end{prop}

Let us introduce now the following notion.
\begin{defn}
The  \textit{Legendre transformation for an optimal control problem} is the mapping
\begin{equation}\label{eq:LegOCP}
{\rm Leg}_{{\mathcal C},{\rm L}}\equiv \tau_{T^*Q}\circ (\alpha_Q^{-1})_{|{\mathcal L}_{{\mathcal C},{\rm L}}} \colon {\mathcal L}_{{\mathcal C},{\rm L}} \longrightarrow T^*Q.
\end{equation}
\end{defn}
It is easy to prove that the optimal control problem is \textit{regular} if and only if ${\rm Leg}_{{\mathcal C},{\rm L}}$ is a local diffeomorphism, see Definition \ref{def:regular:opc}. We will say that the optimal control problem is \textit{hyperregular} if ${\rm Leg}_{{\mathcal C},{\rm L}}$ is a global diffeomorphism.

In the sequel we assume that $\Gamma_{|\widetilde{\mathcal C}}: \widetilde{\mathcal C}\rightarrow {\rm T}Q$ has connected fibers  along the points of its image.
Under this condition we have the following result.
 \begin{lemma}\label{lema-3}
 For any $v_q$ in $\Gamma (\widetilde{\mathcal C})$ 
 \[
 L(c)=L(\tilde{c}) 
 \]
 for all $c$ and $\tilde{c}$ in $\Gamma^{-1}(v_q)\cap \widetilde{\mathcal C}$.
  \proof
For $v_q\in \Gamma(\widetilde{\mathcal C})$ there exists $\tilde{c}$ in $\widetilde{\mathcal C}$ such that $\Gamma(\tilde{c})=v_q$.  Since we assume that  $\Gamma_{|\widetilde{\mathcal C}}: \widetilde{\mathcal C}\longrightarrow TQ$ has connected fibers,
 for any initial point $c$ in  $\Gamma^{-1}(v_q)\cap \widetilde{\mathcal C}$  there exists a curve $\gamma: [0,1]\longrightarrow \widetilde{\mathcal C}$ such that 
 $\gamma(0)=c$, $\gamma(1)=\tilde{c}$ and $\Gamma_{|\widetilde{\mathcal C}}(\gamma(t))=v_q$. Hence
 \[
 \dot{\gamma}(t)\in \ker {\rm T}_{\gamma(t)}\Gamma
 \]
From the definition of $\widetilde{\mathcal C}$ we deduce that 
 \[
 \langle {\rm dL}(\gamma(t)), \dot{\gamma}(t)\rangle=0 \quad \mbox{for all } t\in [0,1].
 \]
 Therefore the Lagrangian $L$ is constant along the curve $\gamma$, that is, $L(\gamma(t))=\hbox{constant}$. Thus, ${\rm L}(c)={\rm L}(\tilde{c})$
 for all $c$ and $\tilde{c}$ in $\Gamma^{-1}(v_q)\cap \widetilde{\mathcal C}$.
   \qed
 \end{lemma}
  
 Now, it will be obtained the ``lagrangian version" of the equations of motion on ${\mathcal L}_{{\mathcal C},{\rm L}}$ for an optimal control problem.
 Any element $v\in {\rm T}_{q}Q$  determines a vertical
vector  at  any point  $w$ in the fiber  over $q$. That vector lying in ${\rm T}_wTQ$ is called the \textit{vertical lift of $v$ at  $w$} and it is denoted
by $v^V_w$. The vertical lift of $v$ at $w$ is the tangent vector
at $t=0$ to the curve $w+t\,v$.
\[
v^V_{w}=\left.\dfrac{\rm d}{{\rm d}t}\right|_{t=0}\left(w+t\,v\right).
\]
Having in mind this definition, it is possible to define now the \textit{energy function} $E_{{\mathcal C},{\rm L}}: {\mathcal L}_{{\mathcal C},{\rm L}}\to \R$ by
\begin{equation}
E_{{\mathcal C},{\rm L}}(\mu_{v_q})=\langle \mu_{v_q}, (v_q)^V_{v_q}\rangle-L(c_q), \quad \mu_{v_q}\in ({\mathcal L}_{{\mathcal C},{\rm L}})_{v_q}, \; c_q\in \widetilde{\mathcal C}, \; \Gamma(c_q)=v_q.
\label{Def:EL}
\end{equation}
Lemma \ref{lema-3} guarantees that the energy is well-defined. 

The following proposition can be easily proved locally.
\begin{prop} Let  $H: T^*Q\times_Q {\mathcal C}\rightarrow \R$ be the Pontryagin's Hamiltonian and $c_q\in \widetilde{\mathcal C}$,
\[
E_{{\mathcal C},{\rm L}}(\mu_{v_q})=H({\rm Leg}_{{\mathcal C},{\rm L}}(\mu_{v_q}), c_q).
\]
\end{prop}

If the mapping $\Gamma: C\rightarrow TQ$ is an embedding, then we can select coordinates in such a way that the control equations are rewritten as follows:
\begin{eqnarray*}
\dot{q}^a&=&u^a,\\
\dot{q}^{\alpha}&=&F^{\alpha}(q^i, u^a).
\end{eqnarray*}
Then the subset ${\mathcal L}_{{\mathcal C},{\rm L}}$ is locally given by
\begin{equation}
\begin{array}{rcl}
a_i+b_{\alpha}\ds{\frac{\partial F^{\alpha}}{\partial q^i}}&=&\ds{\frac{\partial {\rm L}}{\partial q^i}},\\[4mm]
b_a+b_{\alpha} \ds{\frac{\partial F^{\alpha}}{\partial u^a}}&=&\ds{\frac{\partial {\rm L}}{\partial u^a}}.
\end{array}\label{eq:LocalGammaCLExample}
\end{equation}
Note that the submanifold ${\mathcal L}_{{\mathcal C},{\rm L}}$ can be locally described by the coordinates $(q^i,u^a,b_\alpha)$.

 In this case the Pontryagin's Hamiltonian is $H=p_au^a+p_\alpha F^\alpha-L$. Note that the condition for being a Morse family is trivially satisfied because
 \begin{equation*}
 \left(\frac{\partial^2 H}{ \partial p_i\partial u^b}\right)=\begin{pmatrix}
\dfrac{\partial^2 H}{ \partial p_a\partial u^b}& \dfrac{\partial^2 H}{\partial p_{\alpha}\partial u^b} \end{pmatrix}= \begin{pmatrix}
{\rm  I}& \dfrac{\partial \Gamma}{\partial  u^b} \end{pmatrix}
\end{equation*}
has maximum rank.


Due to the local expression of the submanifold ${\mathcal L}_{{\mathcal C},{\rm L}}$, the energy function in~\eqref{Def:EL} is given by
\[E_{{\mathcal C},{\rm L}}(q^i, u^a, b_{\alpha})=u^a\dfrac{\partial L}{\partial u^a}(q^i,u^a)-u^ab_{\alpha}\frac{\partial F^{\alpha}}{\partial u^{a}}(q^i, u^a)+b_{\alpha}F^{\alpha}(q^i, u^a)-L(q^i, u^a).
\]
The Legendre transformation in the coordinates $(q^i,u^a,b_\alpha)$ for  ${\mathcal L}_{{\mathcal C},{\rm L}}$ is written as follows
\[
{\rm Leg}_{{\mathcal C},{\rm L}}(q^i, u^a, b_{\alpha})=\left(q^i,\dfrac{\partial {\rm L}}{\partial u^a}-b_{\alpha} \dfrac{\partial F^{\alpha}}{\partial u^a}, b_{\alpha}\right)
\]
Then the {\bf Poincar\'e-Cartan 2-form} $\Omega_{{\mathcal C},{\rm L}}$ on ${\mathcal L}_{{\mathcal C},{\rm L}}$ is defined by
\[
\Omega_{{\mathcal C},{\rm L}}=({\rm Leg}_{{\mathcal C},{\rm L}})^*\omega_Q
\]
and the ``lagrangian version" of the equations of motion for an optimal control problem is given by the following result.
\begin{prop} 
The Pontryagin's Hamilton equations are equivalent to the following presymplectic system
\[
i_{X}\Omega_{{\mathcal C},{\rm L}}={\rm d} E_{{\mathcal C},{\rm L}}.
\] \label{Prop:LagrangeOCP}
\end{prop}

\subsection{Example: The cheapest stop of a train~\cite[Section 13.3]{Agrachev}}

Let us consider the following control system
\begin{eqnarray}
\dot{x}^1&=&x^2, \label{eq:Hx1}\\ \dot{x}^2&=& u,
\label{eq:Hx2} \end{eqnarray}
which describes the motion of a train with one control. The aim is to stop the train at a fixed final time with a minimum energy. This energy is considered proportional to the integral of the square acceleration, that is, we want to find controlled trajectories that maximize the following functional
\begin{equation*}
\int_0^TL(x^1(t),x^2(t),u(t))\, {\rm d} t=\dfrac{1}{2} \int_0^T u^2(t) \, {\rm d}t. \label{eq:CostTrain}
\end{equation*}
Pontryagin's Hamiltonian~\cite{Pontryagin} is given by
\begin{eqnarray*}
H\colon T^*\mathbb{R}^2\times \mathbb{R}& \longrightarrow &\mathbb{R}, \\
(x,p,u)& \longmapsto & p_1x^2+p_2 u-\dfrac{1}{2} u^2.
\end{eqnarray*}
Hamilton's equation are given by~\eqref{eq:Hx1},~\eqref{eq:Hx2} and
\begin{eqnarray*}
\dot{p}_1&=& 0, \\ \dot{p}_2&=& -p_1.
\end{eqnarray*}
A necessary condition for maximizing the Hamiltonian over the controls is
\begin{equation}
\dfrac{\partial H}{\partial u}=p_2-u=0. \label{eq:HuZero}
\end{equation}
Let us compute the matrix in~\eqref{Eq:SecondDerivativeHu}:
\begin{equation*}
\begin{pmatrix}\ds{\frac{\partial^2 H}{\partial q^i\partial u^b}}  &  \ds{\frac{\partial^2 H}{\partial p_i\partial u^b}}  &  \ds{\frac{\partial^2 H}{\partial u^a\partial u^b}} \end{pmatrix}=\begin{pmatrix} 0 & 0 & 0 & 1 & -1 \end{pmatrix}.
\end{equation*}
Hence, the matrix has maximum rank and the optimal control is regular in the sense of Definition~\ref{Defn:RegOCP}. This example satisfies the typical condition in the literature on the maximum rank of  $\ds{\frac{\partial^2 H}{\partial u^a\partial u^b}}$. In~\cite[Section 13.3]{Agrachev} this problem is solved from there.

However, the above matrix also satisfies that $ \ds{\frac{\partial^2 H}{\partial p_i\partial u^b}}$ has maximum rank. Thus, by implicit function theorem, we can rewrite Hamilton's equations in terms of the local coordinates $(x^1,x^2,p_1,u)$. This local approach arises naturally when we consider the Lagrangian formalism for an optimal control described in the above section. 
Note that 
\begin{equation*}
D_H=\{(x^1,x^2,p_1,u,x^2,u,0,-p_1)\} \quad {\rm and} \quad {\mathcal L}_{{\mathcal C},{\rm L}}=\{(x^1,x^2,x^2,u,0,-p_1,p_1,u)\}.
\end{equation*}
The local coordinates for $ {\mathcal L}_{{\mathcal C},{\rm L}}$ are $(x^1,x^2,p_1,u)$ and the set $\widetilde{\mathcal C}$ in Lemma~\ref{Lemma:TildeC} is equal to ${\mathcal C}$. The condition $\pi_{TQ}({\mathcal L}_{{\mathcal C},{\rm L}})\subseteq \Gamma({\mathcal C})$ in Lemma~\ref{Lemma:TildeC} is
satisfied with equality.

The Legendre transformation for this optimal control problem is locally given by
\begin{equation*}
{\rm Leg}_{{\mathcal C},{\rm L}}(x^1,x^2,p_1,u)=(x^1,x^2,p_1,u),
\end{equation*}
and the energy function is given by 
\begin{equation*}
E_{{\mathcal C},{\rm L}}(x^1,x^2,p_1,u)=H(x^1,x^2,p_1,u,u)=p_1x^ 2+\dfrac{1}{2}u^2.
\end{equation*}
Thus the Poincar\'e-Cartan form $\Omega_{{\mathcal C},{\rm L}}$ on $ {\mathcal L}_{{\mathcal C},{\rm L}}$ is given by
\begin{equation*}
\Omega_{{\mathcal C},{\rm L}}=({\rm Leg}_{{\mathcal C},{\rm L}})^*\omega_Q=-{\rm d}p_1 \wedge {\rm d}x^1-{\rm d}u \wedge {\rm d}x^2.
\end{equation*}
If we write now the ``lagrangian version" of the Hamilton's equations for this optimal control problem given in Proposition~\ref{Prop:LagrangeOCP}, $i_{X}\Omega_{{\mathcal C},{\rm L}}={\rm d} E_{{\mathcal C},{\rm L}}$,
we have
\begin{equation*}
\begin{array}{rclrcl}  \dot{x}^1&=&x^2, &  \dot{p}_1&=& 0, \\ \dot{x}^2&=& u, &\dot{u}&=& -p_1.
\end{array}
\end{equation*}

\subsection{An example of overactuated control system}

Consider now the following control system
\begin{eqnarray*}
\dot{x}&=&u_1, \\
\dot{y}&=& u_2+u_3.
\end{eqnarray*}
The cost function whose functional must be minimized is ${\rm L}(x,y,u_1,u_2,u_3)=\dfrac{1}{2}\left( u_1^2+u_2^2+u_3^2\right)$.

The corresponding Pontryagin's Hamiltonian is $H(x,y,p_1,p_2,u_1,u_2,u_3)=p_1u_1+p_2(u_2+u_3)-\dfrac{1}{2}\left(u_1^2+u_2^2+u_3^2\right)$. In this example the matrix in~\eqref{Eq:SecondDerivativeHu} is given by
\begin{equation*}
\begin{pmatrix}
0 & 0 & 1 & 0 & -1 & 0  & 0 \\
0 & 0 & 0 & 1 & 0 & -1  & 0 \\
0 & 0 & 0 & 1 & 0 & 0  & -1
\end{pmatrix}
\end{equation*}
because
\begin{equation*}
\dfrac{\partial H}{\partial u_1}= p_1-u_1, \quad
\dfrac{\partial H}{\partial u_2}= p_2-u_2, \quad
\dfrac{\partial H}{\partial u_3}= p_2-u_3.
\end{equation*}
 Obviously, the matrix has maximum rank and the optimal control problem is regular because of Definition~\ref{Defn:RegOCP}.
This is an example of overactuated control system, where the controls are not in principle repetitive. However, if we work out the equations defining ${\mathcal L}_{{\mathcal C},{\rm L}}$, we will see that $b_1=u_1$, $b_2=u_2=u_3$.
Note that 
\begin{equation*}
{\mathcal L}_{{\mathcal C},{\rm L}}=\{(x,y,u_1,2u_2,0,0,u_1,u_2)\}.
\end{equation*}
The local coordinates for $ {\mathcal L}_{{\mathcal C},{\rm L}}$ are $(x,y,u_1,u_2)$ and the set $\widetilde{\mathcal C}=\{(x,y,u_1,u_2,u_2)\}\subset {\mathcal C}$. Then Lemma~\ref{Lemma:TildeC} is satisfied with equality
\begin{equation*}
\pi_{TQ}({\mathcal L}_{{\mathcal C},{\rm L}})=\{((x,y,u_1,2u_2)\}=\Gamma(\widetilde{\mathcal C})=\{(x,y,u_1,2u_2)\}.
\end{equation*}
The Legendre transformation for this optimal control problem is locally given by
\begin{equation*}
{\rm Leg}_{{\mathcal C},{\rm L}}(x,y,u_1,u_2)=(x,y,u_1,u_2),
\end{equation*}
and the energy function is given by 
\begin{equation*}
E_{{\mathcal C},{\rm L}}(x,y,u_1,u_2)=H(x,y,u_1,u_2,u_1,u_2,u_3)=u_1^2+2u_2^2- \dfrac{1}{2}\left(u_1^2+2u_2^2\right)=\dfrac{1}{2}u_1^2+u_2^2.
\end{equation*}

Thus the Poincar\'e-Cartan form $\Omega_{{\mathcal C},{\rm L}}$ on $ {\mathcal L}_{{\mathcal C},{\rm L}}$ is given by
\begin{equation*}
\Omega_{{\mathcal C},{\rm L}}=({\rm Leg}_{{\mathcal C},{\rm L}})^*\omega_Q={\rm d}x \wedge {\rm d}u_1+{\rm d}y \wedge {\rm d}u_2.
\end{equation*}
If we write now the ``lagrangian version" of the Hamilton's equations for this optimal control problem given in Proposition~\ref{Prop:LagrangeOCP}, $i_{X}\Omega_{{\mathcal C},{\rm L}}={\rm d} E_{{\mathcal C},{\rm L}}$,
we have
\begin{equation*}
\begin{array}{rclrcl}  \dot{x}&=&u_1, &  \dot{y}&=& 2u_2, \\ \dot{u}_1&=&0, &\dot{u}_2&=& 0.
\end{array}
\end{equation*}

\section{Future work}

As future research line, Morse families might enable a geometric description for the optimal control problems defined on manifolds with corners and/or boundaries. Remember that optimal control problems are considered a generalization of variational calculus~\cite{1997SussmannEt} because in the optimal control problems the controls take values in a closed and bounded set. The Morse families seem to be a useful geometric tool in order to avoid to replace the maximization condition of the Hamiltonian over the controls by the weaker condition $\partial H/\partial u=0$ given in~\eqref{eq:HPartialU0}.

\section*{Acknowledgements}

This work has been partially supported by MICINN (Spain)
MTM2008-00689, MTM2009-13383, MTM2009-08166-E, MTM2010-21186-C02-01 and MTM2010-21186-C02-02; 2009SGR1338 from the Catalan government and SOLSUB200801000238 from the Canary Islands government; the European project IRSES-project “GeoMech-246981” and the ICMAT Severo Ochoa project SEV-2011-0087. D. Iglesias acknowledges MICINN for financial support through
the Ramon y Cajal program

\bibliographystyle{plain}
\bibliography{References}

\end{document}